%% file: LM6_body.tex
\documentclass[a4paper,12pt]{amsart}

\usepackage{amsmath,amsfonts,amsthm,amssymb}

\input{definitions}

\numberwithin{equation}{section}
 \numberwithin{thm}{section}

\begin{document}

\title[Finite jet determination of local CR automorphisms]{Finite jet determination of local CR automorphisms through resolution of degeneracies}
\author{Bernhard Lamel}
\address{Universit\"at Wien, Fakult\"at f\"ur Mathematik, Nordbergstrasse 15, A-1090 Wien, \"Osterreich}
\email{lamelb@member.ams.org}%
\author{Nordine Mir}
\address{Universit\'e de Rouen, Laboratoire de Math\'ematiques Rapha\"el Salem, UMR 6085 CNRS, Avenue de
l'Universit\'e, B.P. 12, 76801 Saint Etienne du Rouvray, France}
\email{Nordine.Mir@univ-rouen.fr}
\thanks{The first author was supported by the FWF, Projekt P17111}

\dedicatory{Dedicated to M.\ Salah Baouendi on the occasion of his
seventieth birthday}

\begin{abstract}
Let $M\subset \C^N$ be a connected real-analytic hypersurface whose
Levi form is nondegenerate at some point. We prove that for every
point $p\in M$, there exists an integer $k=k(M,p)$ such that germs
at $p$ of local real-analytic CR automorphisms of $M$ are uniquely
determined by their $k$-jets (at $p$). To prove this result we develop a
new technique that can be seen as a resolution of the degeneracies of $M$. This procedure consists of blowing 
up $M$ near an arbitrary point $p\in M$ regardless of its minimality or nonminimality;
then, thanks to the blow-up, the original problem can be reduced to an analogous one for a very special class of
nonminimal hypersurfaces for which one may use known techniques to prove the finite jet determination property of its CR automorphisms.
\end{abstract}

\maketitle


\section{Introduction}

This paper is concerned with the finite jet determination problem for germs of
CR automorphisms of real-analytic hypersurfaces in complex space. 
Our main motivation is the following conjecture that essentially goes back
to the recent work of Baouendi, Ebenfelt and Rothschild \cite{BER1}:

\begin{conj}\label{conjecture}
 Let $M\subset \C^{n+1}$ be a connected real-analytic  holomorphically
nondegenerate  hypersurface, $n\geq 1$. Then for every
$p\in M$, there exists a positive integer $k=k(M,p)$ such that germs
at $p$ of local real-analytic CR automorphisms of $M$ are uniquely
determined by their $k$-jets at $p$.
\end{conj}

Let us recall that a (connected) holomorphically nondegenerate real
hypersurface is a real hypersurface for which there is no germ of a nontrivial
holomorphic vector field tangent to an open piece of $M$ (this notion was
introduced by Stanton \cite{S2}). A solution to the above conjecture 
would provide a completely satisfactory local CR version for real-analytic
hypersurfaces of $\C^{n+1}$ of the classical uniqueness theorem of H.\ Cartan
\cite{Hca} stating that holomorphic self-automorphisms of bounded domains in
$\C^{n+1}$ are uniquely determined by their 1-jet at any point of the source
domain. Indeed, holomorphic nondegeneracy  appears to be
the ``natural'' obstruction to 
finite jet determination, the necessity of the condition being 
observed in \cite{BER1}.

Much progress has been made in recent years toward the solution
to the above mentioned conjecture and a number of important cases
have been settled. Historically, the first case considered was when
the given hypersurface has everywhere nondegenerate Levi-form. This was
solved by E. Cartan \cite{Ca1, Ca2}, Tanaka \cite{Ta1} and
Chern-Moser \cite{CM}, as a consequence of their solution to the
biholomorphic equivalence problem. Furthermore, in that setting
unique determination by $2$-jets holds at every point.

In order to follow the lines of the most recent developments related to
Conjecture~\ref{conjecture}, it is relevant to understand the structure of a
real-analytic holomorphically nondegenerate hypersurface and, in particular,
all possible types of Levi-degenerate points that such a manifold can have. And
indeed such a hypersurface may contain points that degenerate in various
possible ways. One first possible situation is given by the case of a
holomorphically nondegenerate real-analytic hypersurface that is everywhere
Levi-degenerate. (This may only happen for $n\geq 2$.) The typical example that
illustrates this situation is the tube in $\C^3$ over the light cone that is
given by the set of smooth points of the real-algebraic variety
$V=\{(z_1,z_2,z_3)\in \C^3: ({\rm Re}\, z_1)^2=({\rm Re}\, z_2)^2+({\rm Re}\,
z_3)^2\}$. For this type of hypersurfaces, neither the results nor the
techniques used in the Levi-nondegenerate case are available.  This led
Baouendi, Huang, Ebenfelt and Rothschild to introduce a finer notion of
nondegeneracy, called {\em finite nondegeneracy} \cite{BHR1,BERbook}. Finite
nondegeneracy is more general than Levi-nondegeneracy (e.g.\ the tube of the
light cone is everywhere finitely nondegenerate); the finite jet determination
problem for the class of finitely nondegenerate hypersurfaces was solved, among
other things, in \cite{BER1}.  This however does not solve Conjecture
\ref{conjecture} since a connected holomorphically nondegenerate real-analytic
hypersurface in general only  satisfies 
the finite nondegeneracy condition at every 
point of a Zariski open subset \cite{BERbook}. Therefore,
the results of \cite{BER1} provide a solution to Conjecture~\ref{conjecture}
for all points $p\in M$ outside a certain proper real-analytic subvariety
$\Sigma$ of $M$.

The above mentioned set $\Sigma\subset M$ consists in a certain
sense  of the most degenerate points that a
hypersurface $M$ as in Conjecture~\ref{conjecture} may have. In
order to deal with this thin set of points, one usually stratifies 
the set $\Sigma$ into a real-analytic subset
$\Sigma_1$ and $\Sigma_2=\Sigma\setminus\Sigma_1$; 
the set $\Sigma_1$ consists of the
points in $\Sigma$ through which there passes a complex
hypersurface of $\C^N$ entirely contained in $M$ (i.e. the set
of nonminimal points), the set $\Sigma_2$ being therefore
the set of minimal points of $\Sigma$. The set $\Sigma_1$ consists
more precisely of a locally finite disjoint union of complex hypersurfaces
of $\C^N$ contained in $M$. The stratification comes from the
fact that points in $\Sigma_1$ and $\Sigma_2$ degenerate in
different manners. Furthermore, and most importantly, there is a
very convenient machinery introduced by Baouendi, Ebenfelt and
Rothschild, the so-called Segre set technique (see e.g.\
\cite{RICM}), that has been extremely useful in the study of
mapping problems but is available only at minimal points. 
Using, among other
tools, this technique, Baouendi, Rothschild and the second author
\cite{BMR1} were able to prove Conjecture~\ref{conjecture} for all
minimal points $p\in M$. This constitutes the most general result
to date toward the solution of the conjecture  for arbitrary $n$.
Hence to complete  the proof of Conjecture~\ref{conjecture}, it
remains to deal with the points lying in the remaining subvariety
$\Sigma_1$.

In the two-dimensional case, Ebenfelt, Zaitsev and the first
author \cite{ELZ1} were recently able to treat this set $\Sigma_1$
successfully, which provided a complete solution to
Conjecture~\ref{conjecture} for $n=1$. (Note that in that case,
the holomorphic nondegeneracy assumption on the hypersurface $M$
is in fact equivalent to its Levi-nonflatness.) In order to prove
the finite jet determination property at every point lying on the
subvariety $\Sigma_1$, the authors of \cite{ELZ1} developed a new
approach, firstly initiated by Ebenfelt \cite{E7} in the study of
the regularity of CR mappings, which, roughly speaking, consists
of reducing the original problem to studying the unique jet
determination of solutions of certain singular systems of differential
equations. This approach strongly contrasts with that used
to deal with the set $\Sigma_2$ of minimal points in e.g.\
\cite{BER5, BMR1} and, therefore, up to now, two different
approaches have been used to study the finite jet determination
problem according to the minimality or nonminimality of the base
point.

In this paper, focusing on the class of {\em generically
Levi-nondegenerate} real-analytic hypersurfaces, we provide a
unified approach to study the finite jet determination problem
and, as a consequence, we give a solution to Conjecture~\ref{conjecture}
for this class.

Our main result is the following.

\begin{thm}\label{t:main}
Let $M\subset \C^{n+1}$ be a connected real-analytic hypersurface
whose Levi form is nondegenerate at some point. Then for every point
$p\in M$, there exists an integer $k=k(M,p)$ such that germs at $p$
of local real-analytic CR automorphisms of $M$ are uniquely
determined by their $k$-jets at $p$.
\end{thm}

In Theorem~\ref{t:main}, even though we assume that the hypersurface 
contains one (and therefore a dense open subset of) Levi-nondegenerate
point(s), there are many Levi-degeneracies that the
hypersurface may have. Indeed, the class of generically
Levi-nondegenerate real-analytic hypersurfaces provides a natural
generalization to  ($n+1$)-dimensional complex euclidean space of
the class of real-analytic Levi-nonflat hypersurfaces in $\C^2$.
Consequently, Theorem~\ref{t:main} may be seen
as a generalization of the mentioned result of
\cite{ELZ1}. On the other hand, this class also contains the
important class of real-analytic hypersurfaces containing no
analytic discs for which the finite jet determination property is
already known to hold in view of the results of \cite{BER5} and
for which the methods of this paper offer a completely different
new proof.

As explained above, the only remaining set $\Sigma_1\subset M$ to be
dealt with in Theorem~\ref{t:main} is that of nonminimal points.
However, our proof applies at {\em every} point of a hypersurface
$M$ satisfying the required conditions regardless of the minimality or nonminimality of the base point.
In that respect, our approach is new and differs from the previous ones.

 Our proof consists of blowing up a hypersurface $M$ satisfying the assumptions
 of Theorem~\ref{t:main} near an arbitrary point $p\in M$. In the preimage of
 the blow-up of $M$, we will obtain a {\em nonminimal} real-analytic
 hypersurface $\widehat M\subset \C^N$ through the origin
 (Proposition~\ref{p:construction}) which has the property that all local
 real-analytic CR automorphisms of $M$ can be lifted to CR automorphisms of
 $\widehat M$ provided their $k$-jet at $p$ coincides with that of the identity
 mapping for $k$ sufficiently large (Proposition~\ref{p:lift}). The finite jet
 determination property of $M$ near $p$ will then be reduced to the same
 property for the constructed nonminimal hypersurface $\widehat M$.
 Furthermore, from our construction, we will get an explicit normal form of
 $\widehat M$ near the origin, which itself follows from an explicit normal
 form for $M$ near $p$ (given in Proposition~\ref{p:normalform}) and the
 explicit form of the blow-up. The obtained explicit normal form of the
 nonminimal hypersurface $\widehat M$ will allow us to use the known techniques
 and results of \cite{E7, ELZ1} to conclude that $\widehat M$ has the desired
 finite jet determination property.

Though the content of this paper is focused on Conjecture~\ref{conjecture}, we
should mention that there has been recently a lot of work on several different
aspects of the finite jet determination of CR maps. In addition to the papers
already mentioned above, we also refer the reader to the papers \cite{Be1, Lo,
Be2, kim, E4, L2, Kowthesis, Travis1, EL1, KZ1, L4, LM2, LM3} and the surveys
\cite{BERbull, Zsurvey, RICM} for more detailed discussions on these various
aspects.

The paper is organized as follows. After settling the notation used throughout
the paper in \S \ref{s:notation}, we introduce in \S \ref{s:spclass} a special
class of (germs of) nonminimal real-analytic hypersurfaces of $\C^{n+1}$ and
prove the finite jet determination of their local CR automorphisms by following
mainly the arguments of \cite{E7, ELZ1}. Such a result will be of fundamental
importance at the end of the paper since  we will show through \S
\ref{s:normalform} -- \S \ref{s:final} that the study of the finite jet
determination problem for local CR automorphisms between real-analytic
hypersurfaces satisfying the conditions of Theorem~\ref{t:main} can be reduced
to that of CR automorphisms of the special class of nonminimal hypersurfaces
defined in \S \ref{s:spclass}. Let us also mention that \S \ref{s:spclass} and
\S \ref{s:normalform}--\S \ref{s:final} are completely independent.

\section{Notation}\label{s:notation}

Throughout the paper, given positive integers $r,k,q$ we denote by $J^r_{0,0}(\C^k,\C^q)$ the jet space of order $r$ of germs of local holomorphic maps $h\colon (\C^k,0) \to (\C^q,0)$; the $r$-jet of $h$ at $0$ is denoted by $j_0^rh$.

For every real-analytic hypersurface $M\subset \C^{n+1}$ and for every point
$p\in M$,  we denote by  $\autMp$ the group of germs at the origin of real-analytic
CR automorphisms of $M$ or, equivalently, of biholomorphisms of $\C^{n+1}$ fixing
$p$ and sending  (the germ of) $M$ into itself. In this paper,
the choice of the reference point $p$ will always be fixed (and we will 
usually have $p$ equal to the origin in $\C^{n+1}$). Recall that  a choice of
local holomorphic coordinates $(z,w)\in \C^n\times \C$, vanishing at $p$, is
called {\em normal} if $M$ can be locally given in these coordinates by an
equation of the form \begin{equation}\label{e:form} w=Q(z,\bar z,\bar w)
\end{equation} for some holomorphic function $Q=Q(z,\chi,\tau)$ defined near
the origin in $\C^{2n+1}$ satisfying
\begin{equation}\label{e:normalcoordinates}
  Q(z,0,\tau)=Q(0,\chi,\tau)=\tau,\quad Q(z,\chi,\bar Q (\chi,z,w))=w,
\end{equation}
where here, and for the remainder of the paper, we denote by
$\bar h$ the series obtained by taking complex conjugates of the coefficients
of the power series $h$. It is well-known that such a choice of normal
coordinates always exists (see e.g.\ \cite{BERbook}). In the remainder of the
paper, we will often use the single notation $(z,w)$ for a choice of normal
coordinates that may change in various propositions and lemmas.

\section{Finite jet determination for a special class of nonminimal hypersurfaces}\label{s:spclass}

The appropriate class of nonminimal hypersurfaces under study in this section is given in the following.

\begin{defin}\label{d:goodnonminimal}
A germ of a real-analytic hypersurface $M\subset \C^{n+1}$ through the origin will be called a {\em good nonminimal} hypersurface if there exists a choice of normal coordinates $(z,w)
\in \C^n \times \C$ such that $M$ is given by an equation of the form \eqref{e:form} with the function $Q$ satisfying
$$Q(z,\chi,\tau) = \tau + \tau^m i \inp{z}{\chi} + \tau^{m+1} \Theta (z,\chi,\tau),$$
where $m$ is some positive integer, $\inp{z}{\chi} = \eps_1 \, z_1\chi_1 + \dots + \eps_n \, z_n \chi_n$, $\eps_j\in \{-1,1\}$ and $\Theta$ is some holomorphic function near $0\in \C^{2n+1}$.
\end{defin}

Let us note that a good nonminimal hypersurface $M$ is of $m$-infinite type (in the sense of \cite{Meylan95, E7}) and that it can also be defined in normal coordinates by a real equation of the form
\begin{equation}\label{e:realform}
 \imag w = \varphi(z,\bar{z},\real w) = (\real w)^m \inp{z}{\bar z}  + O\left((\real w)^{m+1}\right),
\end{equation}
where $\varphi$ is some real-analytic function near the origin in $\R^{2n+1}$. For this and further notions and standard facts about nonminimal (or equivalently, infinite type) hypersurfaces, we refer the reader to the papers \cite{Meylan95, E7, LM4}.

The goal of this section is to provide a proof of the following result:

\begin{prop}\label{p:fjdgood} Let $M\subset \C^{n+1}$ be a germ through the
  origin of a good nonminimal real-analytic hypersurface. Then $(M,0)$ has the
  finite jet determination property: There exists a positive integer $K$
  such that if $H\colon (\C^{n+1},0)\to (\C^{n+1},0)$ is a local biholomorphism
  sending $M$ into itself with the same $K$-jet at the origin as that of the
  identity mapping, then $H$ is the identity. 
\end{prop}

\begin{rem} At this point, we should mention that the proof of
  Proposition~\ref{p:fjdgood} could be achieved by essentially combining
  several arguments and results from the papers \cite{E7, ELZ1}. For the
  reader's convenience, we will give a complete proof of this result
  here, which avoids referring to too many results in the literature.  As in
  \cite{E7}, we shall derive a singular complete system of differential
  equations for the automorphisms of the hypersurfaces considered in
  Proposition~\ref{p:fjdgood}.  Due to the special normal form of the
  hypersurfaces, the construction of this complete system has the advantage to
  be shorter and substantially easier than in \cite{E7}; in particular, we are
  not forced to work on the tangent bundle of the manifolds in this paper. A
  suitable adaptation of the arguments in \cite{ELZ1} then gives the
  desired finite jet determination property for the hypersurfaces under
  consideration.  
\end{rem}

For the proof of Proposition~\ref{p:fjdgood}, we fix a choice of normal
coordinates $(z,w)$ for $(M,0)$ satisfying the conditions of
Definition~\ref{d:goodnonminimal}, and in these coordinates, we split every
$H\in \autM$ as follows $H=(F,G)$. We also use the notation used in the
definition. We start with the following known fact.

\begin{lem}
  \label{l:derstuff} Let $M$, $m$ and $H=(F,G)$ be as given above. Then
  \begin{equation}
    G_{w^\ell} (z,0) = G_{w^\ell} (0) \in \R, \quad \ell \leq m.
    \label{e:gwreal}
  \end{equation}
\end{lem}
\begin{proof} Since  $H$ sends $M$ into itself, we have the identity
  \begin{equation}
    G(z,Q(z,\chi,\tau)) = Q (F(z,Q(z,\chi,\tau)), \bar{F} (\chi,\tau),
    \bar{G} (\chi,\tau)),
    \label{e:basicequ}
  \end{equation}
  which holds for $(z,\chi,\tau)\in\C^{2n+1}$ close enough to $(0,0,0)$.
Since $M$ is a good nonminimal hypersurface and since $H$ is invertible, we have the following relations $$G(z,Q(z,\chi,\tau))=G(z,\tau)+i\tau^m G_w(z,\tau)\inp{z}{\chi}+O(\tau^{m+1})$$ and
$$Q(F(z,Q(z,\chi,\tau)),\bar F(\chi,\tau),\bar G(\chi,\tau))=\bar G(\chi,\tau)+i\bar G(\chi,\tau)^m \inp{F(z,0)}{\bar F(\chi,0)}+O(\tau^{m+1}).$$
 Therefore expanding both sides of \eqref{e:basicequ} as power series in $\tau$, and comparing
  the coefficients of $\tau^\ell$ for $\ell < m$, we get
  \[ G_{w^\ell} (z,0)  = \bar{G}_{w^\ell} (\chi,0); \]
  thus, \eqref{e:gwreal} follows for $\ell < m$. Comparing the
  coefficient of $\tau^m$ in \eqref{e:basicequ}, we also obtain
  \[ G_{w^m} (z,0) + G_w (z,0) i\inp{z}{\chi} =
  \bar{G}_{\tau^{m}} (\chi,0)
  + i\bar{G}_\tau (\chi,0)^m \inp{F(z,0)}{\bar{F}(\chi,0)};\]
  Now setting $\chi = 0$ we obtain \eqref{e:gwreal}
  for $\ell = m$. This completes the proof of the lemma.
\end{proof}

Keeping the notation defined above, for every $H\in \autM$, thanks to Lemma~\ref{l:derstuff} we may write $G(z,w) = P(w) + w^{m} G_2 (z,w)$ where $P(w)=\sum_{j=1}^{m-1}\frac{G_{w^j}(0)}{j!}w^j$. Note that $P(w) = \bar{P} (w)$ and that if $m>1$, $P(w) = w \tilde{P} (w)$ with $\tilde{P} (0)\neq 0$ (for $m=1$, we have $\tilde P \equiv 0$). When $H$ is the identity mapping, the above splitting is written as follows ${\rm Id}=(F^0,G^0)$ with $G^0(z,w)= P^0(w)+w^mG_2^0(z,w)$.  (Note that $G_2^0\equiv 0$ if $m>1$ and $G_2^0\equiv 1$ if $m=1$.)  Observe furthermore that for every $H\in \autM$, we may write $P(\bar w)-P(w)=(\bar w-w)\, {\mathcal Q}(w,\bar w,j_0^{m-1}G)$ and $\tilde P(w)={\mathcal T}(w,j_0^{m-1}G)$ where ${\mathcal Q}, {\mathcal T}$ are universal polynomials in all their arguments.

We will first proceed to derive a singular complete system for the
mapping $V:=(F,G_2)|_M$ associated to any $H\in \autM$ provided
the $(m+2)$-jet of $H$ (at $0$) is sufficiently close to that of
the identity mapping. In the complete system, we will keep the
coefficients of $P$, or equivalently the $(m-1)$-jet of $G$ at
$0$, as parameters. To this end, we need to introduce some further
notation. Since $M$ may also be defined near the origin by an
equation of the form \eqref{e:realform}, we may use
$t:=(z,\bar{z},s)$ where $s=\real w$ as real-analytic coordinates
for $M$ near $0$ and we  also write $\theta := (\real z,\imag z)$.
The precise singular complete system we need is given in the
following lemma.

\begin{lem}\label{l:reflecid}
In the above setting, there exists a real-analytic map $\Upsilon$ defined in a neighbourhood of the point $\varpi_0:=(0,j_0^{m-1}G^0,(((s^m\partial_s^c)\partial_\theta^dV^0)(0))_{c+|d|\leq 2})$
 such that for every $H\in \autM$ with $j_0^{m+2}H$ sufficiently close to $j_0^{m+2}{\rm Id}$, the following identity holds for $(\theta,s)$ sufficiently close to $0\in \R^{2n+1}$:
\begin{equation}\label{e:completeend}
 \left (   \left( s^m \dop{s} \right)^a
 \vardop{}{\theta}{b} V \right)_{a+|b|\leq 3}=
  \Upsilon \left (\theta,s, j_0^{m-1}G,\left( \left(s^{ m}  \dop{s} \right)^c
  \vardop{}{\theta}{d}  V \right)_{c+|d|\leq 2 } \right).
\end{equation}
Here $V=(F,G_2)|_M$ and $V^0:=(F^0,G_2^0)|_M$ are as defined above and $t=(\theta,s)$ are used as local coordinates for  $M$ near $0$.
\end{lem}

\begin{proof}
For every $H\in \autM$, since $H$ sends $M$ into itself, we have the following ``basic equality''
\begin{equation}
\begin{aligned}
  G(z,w) &= Q(F(z,w),\bar{F}(\bar{z},\bar{w}), \bar{G}(\bar{z},\bar{w}))\\
& =
  \bar{G}(\bar z,\bar w) + i \left( \bar{G}(\bar z,\bar w) \right)^{m} \inp{F(z,w)}{\bar F(\bar z,\bar w)}+ \left(\bar{G}(\bar z,\bar w)\right)^{m+1} \Theta (F(z,w),\bar H(\bar z,\bar w)),
\end{aligned}
\label{e:bequ}
\end{equation}
which is valid for  $(z,w)\in M$ close to $0$. Using the local real-analytic coordinates $t$, then \eqref{e:bequ} is
valid for all $t$ in a neighbourhood of $0$ in $\R^{2n+1}$, where we now
think of $z,w,\bar{z},\bar{w}$ as real-analytic functions of $t$. Since $M$  is of
$m$-infinite type, $m\geq 1$, the functions $w=w(t)$ and $\bar{w}=\overline{w(t)}$ have the following properties:
\begin{enumerate}
\item[(P1)] \label{i:dividethrough} The function $A(t):={\overline{w(t)}}/{{w}(t)}$ is  real-analytic (near $0$) and its value at the origin given by $1$;
\item[(P2)] \label{ii:diffdiv} We have $w - \bar{w}= O(s^m)$; this implies that the function $B(t):=(\overline{w(t)} - {w(t)})/{(w(t))^m}$ is also real-analytic near $0$.
\end{enumerate}

We rewrite the basic equation \eqref{e:bequ} as
\begin{multline}\label{e:bequ2}
w^{m} G_2(z,w) = \left( {P} (\bar{w}) - {P}(w) \right) + \bar{w}^{m} \bar{G}_2(\bar z,\bar w) + i \left( \bar{w}^m \right)\left( \bar{\tilde{P}}(\bar w) +
 \bar{w}^{m-1} \bar{G}_{2}(\bar z,\bar w) \right)^m \inp{F(z,w)}{\bar F(\bar z,\bar w)}\\
 +\bar{w}^{m+1} \left( \bar{\tilde{P}}(\bar w) +
 \bar{w}^{m-1} \bar{G}_{2}(\bar z,\bar w) \right)^{m+1}
\Theta (F(z,w),\bar F(\bar z,\bar w), \bar w \bar{\tilde{P}}(\bar w)+\bar w^m \bar{G}_2(\bar z,\bar w)),
\end{multline}
which holds for all $t\in \R^{2n+1}$ sufficiently close to $0$. In what follows,  if the variables are not
written, it is understood that barred functions have $(\bar{z}, \bar{w})$ as their arguments, while
unbarred functions  have $(z,w)$ as arguments. It is also understood that we still view the functions $(z,\bar z,w,\bar w)$ as functions of $t$, but we do not write it in order to avoid too heavy equations.

In view of Properties (P1) and (P2), we may divide \eqref{e:bequ2} by $w^m$ and obtain an identity of the form
\begin{multline}
  G_2 = B(t)\, {\mathcal Q}(w,\bar w,j_0^{m-1}G) + \left( A(t)\right)^m \bar{G}_2 +
  i \left( A(t)\right)^m
  \left( \bar{{\mathcal T}}(\bar w,j_0^{m-1}G) + \bar{w}^{m-1}\bar{G}_2 \right)^m \inp{F}{\bar{F}} +\\
  \bar{w}  \left( A(t) \right)^m  \left( \bar{{\mathcal T}}(\bar w,j_0^{m-1}G) +
 \bar{w}^{m-1} \bar{G}_{2} \right)^{m+1}
\Theta (F,\bar F, \bar w \bar{{\mathcal T}}(\bar w,j_0^{m-1}G)+\bar w^m \bar{G}_2),
  \label{e:bequ23}
\end{multline}
valid for all $t\in \R^{2n+1}$ sufficiently close to $0$. Notice that we may rewrite the expression
$$\left( \bar{{\mathcal T}}(\bar w,j_0^{m-1}G) +
 \bar{w}^{m-1} \bar{G}_{2} \right)^{m+1}
\Theta (F,\bar F, \bar w \bar{{\mathcal T}}(\bar w,j_0^{m-1}G)+\bar w^m \bar{G}_2)$$
in the following form
$$R(t,F,\bar F, \bar G_2,j_0^{m-1}G)$$
for some universal function
$R$ in all its arguments (depending only on $M$) holomorphic in some neighbhourhood of $\{0\}\times \C \times J_{0,0}^{m-1}(\C^{n+1},\C)\subset \C^{4n+1} \times \C \times J_{0,0}^{m-1}(\C^{n+1},\C)$.

We now proceed by differentiating the identity
\begin{multline}\label{e:bequ3}
  G_2 = B(t)\, {\mathcal Q}(w,\bar w,j_0^{m-1}G) + \left( A(t)\right)^m \bar{G}_2 +\\
  i \left( A(t)\right)^m
  \left( \bar{{\mathcal T}}(\bar w,j_0^{m-1}G) + \bar{w}^{m-1}\bar{G}_2 \right)^m \inp{F}{\bar{F}} +
  \bar{w}  \left( A(t) \right)^m  R(t,F,\bar F,\bar G_2,j_0^{m-1}G)
\end{multline}
with respect to the CR vector
fields on $M$. To this end, with our choice of coordinates for $M$, we choose as a basis of the CR structure on $M$ the following vector fields
\[ L_j :=\dop{\bar{z}_j}
- \frac{ \varphi_{\bar{z}_j}}{\varphi_s - i}\, \dop{s},\quad j=1,\ldots,n. \]
In what follows, a multiindex $J$ is a finite sequence $(J_1,\ldots,J_\ell)\subset \{1,\ldots,n\}^\ell$ and we denote by $|J|:=\ell$ the length of $J$. Then $L^{J}$ (resp.\ $\bar L^{J}$) denotes the differential operator $L_{J_1}\ldots L_{J_\ell}$ (resp.\ $\bar{L}_{J_1}\ldots \bar{L}_{J_\ell}$).

We now apply $L_j$ to \eqref{e:bequ3} for every $j=1,\ldots,n$. Since $L_j G_2 = 0$, we obtain an equation of the form
\begin{equation}\label{e:system}
 0 = \Psi_j \left(t,F,\bar{F}, \bar G_2, L_j \bar F, L_j \bar{G}_2,j_0^{m-1}G \right),
\end{equation}
for $t\in \R^{2n+1}$ close to $0$, where $\Psi_j$ is holomorphic in a neighbourhood of
$\{0\}\times \C \times \C^n \times \C \times J_{0,0}^{m-1}(\C^{n+1},\C)\subset \C^{6n+3}\times J_{0,0}^{m-1}(\C^{n+1},\C)$ (and independent of the mapping $H$). We denote by $A=(A_1,\ldots,A_n)$ the second set of variables appearing in the function $\Psi_j$ (and in which the map $F$ appears). We now evaluate the partial derivative of $\Psi_j$ with respect to $A_k$ at the reference point corresponding to the identity mapping. Using the fact that $L_j \bar{w} (0) = 0$ and $A(0)=1$,  we obtain
\[ \dopt{\Psi_j}{A_k} (0,0,0,\bar G_2^0(0),L_j \bar{F}^0 (0),L_j \bar{G}_2^0(0),j_0^{m-1}G^0) =
i \eps_k \delta_{j,k}, \]
where $\delta_{j,k}$ denotes the Kronecker symbol. Thus, we may apply 
the Implicit Function Theorem to the equations \eqref{e:system} for
$j=1,\ldots,n$ and conclude that there exists
a $\C^n$-valued holomorphic map $\Phi$ defined in a neighbourhood of $\eta_0:=(0,\left((L^\alpha \bar F^0 (0),L^\alpha \bar G_2^0(0))\right)_{|\alpha|\leq 1},j_0^{m-1}G^0)$ in $\C^{2n+1}\times \C^{(n+1)^2}\times J_{0,0}^{m-1}(\C^{n+1},\C)$ such that for every $H\in \autM$ whose $(m+1)$-jet at $0$ is sufficiently close to that of the identity mapping, then
\begin{equation}\label{e:el}
 F = \Phi \left(t,\left(L^{\alpha} \bar{F}, L^{\alpha} \bar{G}_2\right)_{|\alpha|\leq 1},j_0^{m-1}G \right), \end{equation}
for $t\in \R^{2n+1}$ close to $0$.
 Plugging \eqref{e:el} into \eqref{e:bequ3} we get an analogous equation for $G_2$. For every map $H\in \autM$ as above, recalling that  $V = (F, G_2)|_M$, we get that $V$ satisfies an identity in the neighbourhood of the origin in $\R^{2n+1}$ of the form
\begin{equation}
  \label{e:startcomplete} V = \Xi \left(t,\left(L^{\alpha} \bar{V} \right)_{|\alpha|\leq 1},j_0^{m-1}G \right),
\end{equation}
where $\Xi$ is a $\C^{n+1}$-valued holomorphic map defined in a neighbourhood of $\eta_0$. In what follows, the elements  $H\in \autM$ for which $j_0^{m+1}H$ is close enough to $j_0^{m+1}{\rm Id}$ so that \eqref{e:startcomplete} holds will be called admissible maps.

Having established \eqref{e:startcomplete}, we can   now use arguments similar to
\cite{E7} but somewhat easier and closer to the paper of Han \cite{HAN}
in order to derive the desired singular
complete system. For this, we set
\[ S := s^m \dop{s},\]
and note that we have the commutation identities
\begin{equation}\label{e:commutation}
 \left[ L_j, L_k \right] = 0, \quad \left[ L_j, \bar{L}_k \right]= a_{j,k} S,
\quad \left[ L_j,S \right] = b_j S,\quad j,k=1,\ldots,n,
\end{equation}
where the $a_{j,k}$ and the $b_j$ are germs of real-analytic functions
at $0\in \R^{2n+1}$, and $a_{j,j} (0) \neq 0$.

Let us now recall the following fact
from \cite[Proposition 4.3]{E7}: for any multi-index $J$, integer $k\geq 1$, there exist real-analytic functions  $b^{e_1\ldots
e_m}_q$ near $0$ in $\R^{2n+1}$ such that
\begin{equation}\label{e:maincommutatoridentity}
\sum_{m=1}^{|J|+k}\sum_{q=0}^{k} b_q^{e_1\ldots e_m}
\underbrace{[\ldots [\bar{L}_{e_1}\ldots \bar{L}_{e_m},L_{1}],L_{1}]\ldots,L_{1}]}_{\text{\rm length $q$ }}= \bar{L}^{J}S^k.
\end{equation}
Here the length of the commutator $[\ldots [X,Y_{1}],Y_{2}]\ldots,Y_{q}]$ is $q$.

 For every integer $j$ and multiindex $J$, we use the notation
 $\Lambda_{j,J}$ to denote new variables lying in $\C^{n+1}$ and $\lambda$
 for a jet in the space $J_{0,0}^{m-1}(\C^{n+1},\C)$. Further, for every
 nonnegative integers $p$ and $q$, we denote by $C_{p,q}$ (resp.\ $\bar{C}_{p,q}$) the set of
 $\C^{n+1}$-valued holomorphic maps ${\mathcal A}={\mathcal
 A}\left(t,(\Lambda_{j,J})_{j+|J|\leq p \atop j\leq q}, \lambda
 \right)$ that are polynomials in the variables $({\Lambda}_{j,J})_{j>0}$
 and $(\Lambda_{0,J})_{|J|>1}$  with holomorphic coefficients in a
 neighbourhood of  $\omega_0:=(0,\left((L^J \bar F^0 (0),L^J \bar
 G_2^0(0))\right)_{|J|\leq 1},j_0^{m-1}G^0)$ in $\C^{2n+1}\times
 \C^{(n+1)^2}\times J_{0,0}^{m-1}(\C^{n+1},\C)$ (resp.\  in a neighbourhood of $\bar{\omega}_0$).

 Applying $\bar{L}_j$ to \eqref{e:startcomplete}
and  using the fact that $\bar{L}_j\bar V=0$,
we get that for every $j=1,\ldots,n$, there exists a map $\Xi_j\in C_{1,1}$ such
 for all admissible maps $H$
\begin{equation}\label{e:bibe}
\bar{L}_j V=\Xi_j(t,(L^J S \bar V)_{j+|J|\leq 1 \atop j\leq 1},j_0^{m-1}G).
\end{equation}
Due to the commutation identities \eqref{e:commutation}, further applications of any vector field $\bar{L}_k$ to \eqref{e:bibe}
will not increase the order of differentiation of the arguments of the right hand side of \eqref{e:bibe},
and therefore for every multiindex $E$, there exists $\Xi_E \in C_{1,1}$, such that for every admissible map $H$, the following identity holds in a neighbourhood of the origin in $\R^{2n+1}$
\begin{equation}\label{e:lastequation}
\bar{L}^{E} V=\Xi_{E}(t,(L^J S^j \bar V)_{j+|J|\leq 1 \atop j\leq 1},j_0^{m-1}G).
\end{equation}
In order to get a similar identity for the map $L^J S^kV$ for every multiindex $J$ and integer $k$, we note that the commutator identity \eqref{e:maincommutatoridentity} and the fact that $L_1 V = 0$ imply
 that we just need to apply
 at most $|J|+k$ vector fields of the family $\{ \bar{L}_1,\ldots, \bar{L}_n \}$ to \eqref{e:startcomplete}, followed
 by an application of at most $k$ instances of $L_{1}$. Therefore to get the desired equations
for the expressions $L^{J} S^{k}V$ we only need to apply at most $k$ instances of the vector field $L_1$ to \eqref{e:lastequation}. We thus get for every multiindex $E$
the existence of $\Xi_{E,k}\in C_{1+k,1}$ such that
\[ L_1^{k} \bar{L}^{E} V = \Xi_{E,k} \left((t,(L^J S^j \bar V)_{j+|J|\leq 1+k \atop j\leq 1},j_0^{m-1}G  \right),\]
for $t$ close to $0\in \R^{2n+1}$ and hence for every
$J$ and $k$, there exists $\tilde{\Xi}_{J,k} \in C_{1+k,1}$ such that
for all admissible maps $H$,
\begin{equation}\label{e:endcomplete}
  \bar{L}^{J} S^{k} V = \tilde{\Xi}_{J,k} \left((t,(L^K S^\ell \bar V)_{\ell+|K|\leq 1+k \atop \ell \leq 1},j_0^{m-1}G  \right).
\end{equation}
Taking the complex conjugate of \eqref{e:endcomplete} for $k = 1$ and plugging the resulting equation into \eqref{e:endcomplete} for
$|J|+ k \leq 3$, we get for any such $J$ and $k$ a map $\hat{\Xi}_{J,k}\in \bar{C}_{2,1}$ such that for all admissible
maps $H$,
\begin{equation}\label{e:endcomplete2}
  \bar{L}^{J} S^{k} V = \hat{\Xi}_{J,k} \left(t,(\bar{L}^K S^\ell  V)_{\ell+|K|\leq 2 \atop \ell \leq 1},j_0^{m-1}G  \right),
\end{equation}
in a neighbourhood of the origin in $\R^{2n+1}$. Since $L_\nu V = 0$ for all $\nu =1,\ldots,n$, rewriting the system \eqref{e:endcomplete2} in the real coordinates $t=(\theta,s)$ we get that each map $V$ as above satisfies a singular complete system of order $3$,
i.e.\  that \eqref{e:completeend} is satisfied for every map $V$ corresponding to a map $H\in \autM$ for which $j_0^{m+2}H$ is sufficiently close to $j_0^{m+2}{\rm Id}$. This completes the proof of Lemma~\ref{l:reflecid}.
\end{proof}

\begin{proof}[Proof of Proposition~{\rm \ref{p:fjdgood}}] Keeping the previously defined notation, we shall prove that there exists an integer $K\geq (m+4)$ such that if $H\in \autM$ satisfies $j_0^{K}H=j_0^K{\rm Id}$ then the complete system satisfied by the corresponding map $V=(F,G_2)|_M$ (and given in Lemma~\ref{l:reflecid}) has only $V^0$ as a solution, which shows that $H={\rm Id}$.

We start by considering only those $H\in \autM$ for which
$j_0^{m+4}H=j_0^{m+4}{\rm Id}$ and making the following
observation. In view of  Lemma~\ref{l:reflecid} all such maps have
their component $V=(F,G_2)|_M$ that satisfy the same complete
system.  For all such $V$'s, the map $V(\theta,0)$ satifies an
ordinary complete system of order $3$ (in the variable $\theta$),
and, therefore, $V(\theta,0)$ is real-analytically parametrized by
its $2$-jet at $0$ (see e.g.\ \cite[Proposition 3.54]{E4}).
Similarly, after differentiating \eqref{e:completeend} with
respect to $s$ $\ell$ times, we see that for every integer $\ell$
the maps $V_{s^\ell}(\theta,0)$ satisfy analogously  an ordinary
complete system of order $3$ and, hence, such maps
$V_{s^\ell}(\theta,0)$ are also real-analytically parametrized by
the $(\ell+2)$-jets of $V$ at $0$.

Next, writing for every $V$ as above
\[ U (\theta,s) :=    \left(\left( s^m \dop{s} \right)^a
  \vardop{}{\theta}{b}  V \right)_{a+|b| \leq 2 },\quad {\mathcal U}(\theta,s):=U(\theta,s)-U(\theta,0),\]
and $U^0$, ${\mathcal U}^0$ for the corresponding maps obtained from $V^0$, we note that it follows from the  observation and the assumption $j_0^{m+4}H=j_0^{m+4}{\rm Id}$ that $U(\theta,0)=U^0(\theta,0)$. It therefore follows from
\eqref{e:completeend} that ${\mathcal U}$ fulfills a real-analytic equation of the form
\[ s^m \frac{\partial {\mathcal U}}{\partial s}(\theta,s) = \widetilde \Upsilon \left(\theta,s, j_0^{m-1}G^0,{\mathcal U}(\theta,s)+U^0(\theta,0) \right), \]
where $\widetilde \Upsilon$ is real-analytic in a neighbhourhood
of $\varpi_0$ (where $\varpi_0$ is given in
Lemma~\ref{l:reflecid}). Since ${\mathcal U}(\theta,0)=0$, we are
now in a position to apply the determination theorem~\cite[Theorem
5.1]{ELZ1} which provides an integer $k$ such that if ${\mathcal
U}$ satisfies ${\mathcal U}_{s^j} (\theta,0) = {\mathcal
U}^0_{s^j} (\theta,0)$ for $j\leq k$, then ${\mathcal
U}(\theta,s)={\mathcal U}^0(\theta,s)$. But, again, in view of the
observation, the last condition is fulfilled provided the
$(k+4)$-jet of $V$ at $0$ agrees with that of $V^0$. Now setting
$K:=m+4+k$, we get the desired result. The proof of the
proposition is complete.
\end{proof}

\section{A normal form for a generically Levi-nondegenerate real-analytic hypersurface}\label{s:normalform}

The goal of this section is to provide, for any connected
real-analytic hypersurface $M\subset \C^{n+1}$ which is
Levi-nondegenerate at some point, a normal form for the hypersurface
near {\em any} of its points. This normal form will be provided by
Proposition~\ref{p:normalform} below.

In the following lemma, the notation $T_p^cM$ denotes the usual
complex tangent space of $M$ at $p$ (see e.g.\ \cite{BERbook}). The
content of the following lemma is well-known, but for the
reader's convenience we provide a proof.

\begin{lem}\label{l:transversecurve}
 Let $M\subset\C^{n+1}$ be a real-analytic hypersurface and $p\in M$. Let $\Gamma$ be a real-analytic curve passing through $p$ that is transverse to $T_p^cM$. Then
 there exists a choice of normal coordinates $(z,w)\in \C^n \times \C$ for $M$ near $p$, vanishing at $p$, such that the germ at $p$ of $\Gamma$ is given in
 these coordinates by the germ of the real line $\{(0,s):s\in (\R,0)\}$.
\end{lem}

\begin{proof}
 We first choose an arbitrary set of normal coordinates $(z',w')\in \C^n \times \C$, vanishing at $p$,
 so that $M$ is given near the origin by an equation of the form
\begin{equation}\label{e:start}
{\rm Im}\, w'=\psi (z',\bar z',{\rm Re}\, w'),
\end{equation}
with $\psi$ real-analytic in a neighborhood of $0\in \R^{2n+1}$
and satisfying $\psi (z',0,{\rm Re}\, w')=\psi (0,\bar z',{\rm
Re}\, w')\equiv 0$. Let $ (\R,0)\ni t\mapsto \gamma (t)\in
(\C^{n+1},0) $ be a parametrization of the curve $\Gamma$ near the
origin. Since $\Gamma$ is transverse to $T_p^cM$ by assumption, we
may assume, after reparametrizing the curve if necessary that
$\gamma (t)=(\beta (t),\eta (t))\in \C^n \times \C$ with $\eta
(t)=t+i\psi (\beta (t),\overline{\beta (t)},t)$ for $t\in \R$
close enough to the origin. In what follows, we shall complexify
the maps $\beta$ and $\eta$ and keep the same notation for the
complexified maps. By the implicit function theorem,
\eqref{e:start} is equivalent to an equation of the form
$w'=Q'(z',\bar z',\bar w')$ for some holomorphic function $Q'$
near the origin in $\C^{2n+1}$ satisfying
\eqref{e:normalcoordinates} (with $Q$ replaced by $Q'$). Consider
the coordinates $(z,w)\in \C^n \times \C$ defined by the following
holomorphic change of variables:
\begin{equation}
 \begin{cases}
  z'=z+\beta (w),\cr
  w'=Q'(z+\beta (w),\bar \beta (w),\bar \eta (w)).
 \end{cases}
\end{equation}
By using the normality of the $(z',w')$ coordinates and the fact
that $\eta (w)=Q'(\beta (w),\bar \beta (w),\bar \eta (w))$ for $w\in
\C$ close to $0$ (which comes from the fact that $\Gamma \subset
M$), the reader can verify that the obtained $(z,w)$ coordinates
satisfy all the requirements of the lemma. The proof of
Lemma~\ref{l:transversecurve} is complete.
\end{proof}

We now go through the construction of the desired normal form for a real-analytic generically Levi-nondegenerate
hypersurface. In what follows, for any matrix $U$ with entries in $\C$, we set $U^*:=\overline{^tU}$, where $^tU$ denotes the transpose matrix of $U$. We will make use of the following result that goes back to Rellich \cite{Rellich} (see also \cite{GinH}):

\begin{prop}\label{p:diagonalization}
Let $A(s)$ be a complex-valued $k\times k$ matrix whose entries depend
real-analytically on $s\in
(a,b)\subset \R$. Assume that  $A(s)$ is hermitian for every $s\in (a,b)$.
Then there exists a $k\times k$ unitary matrix $U(s)$, depending real-analytically on $s\in (a,b)$ such that
the matrix $$U(s)A(s)U^*(s)$$ is a diagonal matrix $($with real-analytic entries in $s)$.
\end{prop}

Proposition~\ref{p:diagonalization} is the main ingredient 
for the construction of the following normal form.

\begin{prop}\label{p:normalform}
Let $M\subset \C^{n+1}$ be a connected real-analytic hypersurface whose Levi form is nondegenerate at some point.
Then, for every $p\in M$, there exist normal coordinates $(z,w)\in \C^n\times \C$, vanishing at $p$, such that $M$ near $p$ is given in these coordinates by an equation of the following form
\begin{equation}\label{e:normalformequation}
 {\rm Im}\, w=\sum_{j=1}^n\eps_j|z_j|^2\, ({\rm  Re}\, w)^{b_j}\, \theta_j({\rm Re}\, w)+R(z,\bar z,{\rm Re}\, w),
\end{equation}
where each $\eps_j\in \{-1,1\}$, $b_j\in \N$, $b_1\geq b_2\geq \ldots \geq b_n$, $\theta_j$ is a real-analytic function near $0\in \R$ satisfying $\theta_j(0)=1$, and $R$ is a real-analytic function near $0\in \R^{2n+1}$ satisfying $$R(z,0,{\rm Re}\, w)=R(0,\bar z,{\rm Re}\, w)\equiv 0,\quad {\rm and}\quad  R(z,\bar z,{\rm Re}\, w)=O(|z|^3).$$
\end{prop}

\begin{proof}
Since the set of Levi-nondegenerate points is dense in $M$, we may
choose a real-analytic curve $\Gamma$ passing through $p$,
transverse to $T_p^cM$, and passing through a Levi-nondegenerate
point of $M$. Shrinking furthermore the curve $\Gamma$ near $p$ if
necessary, we may assume that $p$ is the only point of $\Gamma$
where $M$ may be Levi-degenerate. By Lemma~\ref{l:transversecurve},
we may choose normal coordinates, that we denote by $(z',w')$,
vanishing at $p$, such that $M$ is given by an equation of the form
$${\rm Im}\, w'=\phi (z',\bar z',{\rm Re}\, w'),$$
with $\phi$ real-analytic near the origin in $\R^{2n+1}$, $\phi
(z',0,{\rm Re}\, w')=\phi (0,\bar z',{\rm Re}\, w')\equiv 0$ and
such that the germ at $p$ of the curve $\Gamma$ is given by (the
germ at $0$ of) the real line $\{(0,s):s\in \R\}$. From the
construction, we also know that for $s\in \R$ sufficiently small,
the $n\times n$ matrix $A(s):=\phi_{z\bar{z}}(0,s)$ is nondegenerate
for all $s\not =0$. Since $A(s)$ is hermitian for each $s$ close to
$0$, we may apply Proposition~\ref{p:diagonalization} to get the
existence of a unitary matrix $U(s)$ real-analytic near the origin
satisfying
\begin{equation}
U(s)A(s)U^*(s)={\mathcal D}(s)=\left(\begin{array}{ccc}
D_1(s)&0&0\\
0&\ddots&0\\
0&0&D_n(s)
\end{array} \right),
\end{equation}
with each $D_j$, $j=1,\ldots,n$, being a nonzero real-analytic
function near the origin, with the origin as the only possible zero
for each $D_j$. We now construct the desired normal coordinate
system $(z,w)$ by defining it through the following holomorphic
change of variables:
\begin{equation}
\begin{cases}
 z'=U(w)\cdot z\cr
 w'=w
\end{cases}
\end{equation}
Here again, we also use the notation $U$ to mean the
complexification of $U$ to a sufficiently neighbourhood of $0$ in
$\C$. In the new $(z,w)$ coordinates, the germ of  $M$ at $0$ is
given by the equation
$${\rm Im}\, w=\phi (U(w)\cdot z,\overline{U(w)\cdot z},{\rm Re}\, w),$$
or equivalently, by using the implicit function theorem by another equation of the form
$${\rm Im}\, w=\Phi (z,\bar{z},{\rm Re}\, w),$$
for some other real-analytic function $\Phi$ defined near the
origin in $\R^{2n+1}$. It is easy to check that $(z,w)$ are normal
coordinates for $M$ i.e.\ that $\Phi (z,0,{\rm Re}\, w)=\Phi
(0,\bar{z},{\rm Re}\, w)\equiv 0$ and that the $n\times n$ matrix
$\Phi_{z\bar{z}}(0,s)$ is equal to ${\mathcal D}(s)$. From this,
it is not difficult to derive the final normal form given by
\eqref{e:normalformequation}, rescaling and interchanging the
$z$-coordinates if necessary. The proof of
Proposition~\ref{p:normalform} is therefore complete.
\end{proof}

\section{Resolution of degeneracies via blow-ups and proof of Theorem~\ref{t:main}}\label{s:final}

\subsection{The blow-up procedure}\label{s:blow-up}
Throughout \S \ref{s:blow-up}, we fix a germ of a real-analytic
hypersurface $M$ through the origin in $\C^{n+1}$
 for which there exists a system of normal
coordinates $(z,w)=(z_1,\ldots,z_n,w)\in \C^n\times \C$, that we
also fix for the remainder of this section, such that $M$ is given by an equation  of the form \eqref{e:normalformequation} with all corresponding quantities satisfying the conditions of Proposition~\ref{p:normalform}.

Our goal here is to construct a good nonminimal real-analytic
hypersurface $\widehat M$ (as defined in \S \ref{s:spclass})
contained in the pre-image of $M$ under a certain explicit
blow-up. It should be mentioned that the use of blow-ups in $\C^2$
already appears for other purposes in the papers
\cite{Zsurvey,Kowthesis,LM4}.

The relevant statement for this paper is given by the following.

\begin{prop}\label{p:construction}
Let $M$ be given as above and ${\mathcal B}\colon \C^{n+1}\to
\C^{n+1}$ be the polynomial map given by
\begin{equation}\label{e:defineblowup}
 {\mathcal B}(z,w):= (z_1w^{\alpha_1},\ldots,z_nw^{\alpha_n},w^2),
\end{equation}
with $\alpha_j:=2+3b_1-b_j$ for $j=1,\ldots,n$. Then there exists a
unique germ of a real-analytic hypersurface $\widehat M$ through the
origin in $\C^{n+1}$, contained in ${\mathcal B}^{-1}(M)$,  which
is of the form \begin{equation}\label{e:form1}
{\rm Im}\, w=({\rm Re}\, w)^{3+6b_1}\eta (z,\bar
z,{\rm Re}\, w),
\end{equation}
 for some real-analytic function $\eta$ defined near $0\in \R^{2n+1}$. Furthermore,
$\eta$ satisfies
\begin{equation}\label{e:form2}
\eta (z,\bar z,0)=\frac{1}{2}\sum_{j=1}^n\eps_j|z_j|^2,\quad {\rm and}\quad  \eta (z,0,{\rm Re}\, w)=\eta (0,\bar z,{\rm Re}\, w)=0.
\end{equation}
\end{prop}

\begin{proof} In what follows, we set $s:={\rm Re}\, w$ and $t={\rm Im}\, w$. First note that the set of points $(z,s+it)$ close to the origin that are contained in the preimage of $(M,0)$ under ${\mathcal B}$ is given by the equation
\begin{equation}\label{e:preimagedef}
2st=\sum_{j=1}^n\eps_j|z_j|^2\, (s^2-t^2)^{b_j}\, (s^2+t^2)^{\alpha_j}\, \theta_j (s^2-t^2)+\widetilde R(z,\bar z,s,t),
\end{equation}
 with \begin{equation}\label{e:orderequation}
\widetilde R(z,\bar z,s,t)=R(z_1w^{\alpha_1},\ldots,\bar z_n\bar w^{\alpha_n}, {\rm Re}\, (w^2)).
\end{equation}
Setting  $t=s^{3+6b_1}v$, we will show that there exists a unique real-analytic function $v=v(z,\bar z, s)$ such that $t=s^{3+6b_1}v$ satisfies \eqref{e:preimagedef}. Plugging $t=s^{3+6b_1}v$ in \eqref{e:preimagedef}, we get the following equation
\begin{equation}\label{e:use}
 2s^{4+6b_1}v=\sum_{j=1}^n\eps_j |z_j|^2\, s^{2b_j+2\alpha_j}(1-s^{4+12b_1}v^2)^{b_j}\, (1+s^{4+12b_1}v^2)^{\alpha_j}\, \theta_j (s^2-s^{6+12b_1}v^2)+\widetilde R(z,\bar z,s,s^{3+6b_1}v).
\end{equation}
By our choice of $\alpha_j$, we have
$s^{2b_j+2\alpha_j}=s^{4+6b_1}$. To show that there exists a
unique real-analytic function $v=v(z,\bar z,s)$ satisfying
\eqref{e:use}, it is enough to see that $\widetilde R(z,\bar
z,s,s^{3+6b_1}v)=O(s^{5+6b_1})$. We claim that one has even
$\widetilde R(z,\bar z,s,s^{3+6b_1}v)=O(s^{6+6b_1})$. Indeed,
since $R(z,\bar z,{\rm Re}\, w)=O(|z|^3)$, we have in view of
\eqref{e:orderequation} that $\widetilde R (z,\bar
z,s,t)=O(|w|^{3\alpha_1})=O(s^{3\alpha_1})=O(s^{6+6b_1})$, proving
the claim. Dividing \eqref{e:use} by $s^{4+6b_1}$, we get that $v$
satisfies the real-analytic equation
\begin{equation}\label{e:hire}
2v=\sum_{j=1}^n\eps_j |z_j|^2\, (1-s^{4+12b_1}v^2)^{b_j}\, (1+s^{4+12b_1}v^2)^{\alpha_j}\, \theta_j (s^2-s^{6+12b_1}v^2)+S(z,\bar z,s,v),
\end{equation}
with $S(z,\bar z,s,v)=O(s^2)$. By the implicit function theorem, \eqref{e:hire} has a unique real-analytic function solution $v=v(z,\bar z,s)$. This proves the first part of the proposition.

For the second part of the proposition, we first notice that
$S(z,0,s,v)= \displaystyle \frac{\widetilde
R(z,0,s,s^{3+6b_1}v)}{s^{4+6b_1}}=0$ and similarly $S(0,\bar
z,s,v)=0$. This proves that $\eta (z,0,s)=\eta (0,\bar z,s)=0$.
Finally, the fact that $\eta (z,\bar
z,0)=\frac{1}{2}\sum_{j=1}^n\eps_j |z_j|^2$ clearly follows from
\eqref{e:hire} and the fact that $\theta_\nu (0)=1$ for all
$\nu=1,\ldots,n$. The proof of Proposition~\ref{p:construction} is
complete.
\end{proof}


We now come to the crucial property satisfied by the hypersurface
$\widehat M$ constructed in Proposition~\ref{p:construction}.

\begin{prop}\label{p:lift}
Let ${\mathcal B}$, $M$ and $\widehat M$ be as above. Then for
every integer $\ell\geq {\rm max}\, \{\alpha_n,3+6b_1\}$, where
$\alpha_n$ and $b_1$ are given in Proposition~{\rm
\ref{p:construction}}, if $H\colon (\C^{n+1},0)\to (\C^{n+1},0)$
is a holomorphic mapping sending $M$ into itself satisfying
$j_0^{\ell}H=j_0^{\ell}{\rm Id}$, then there exists a unique
holomorphic map $\widehat H\colon (\C^{n+1},0)\to (\C^{n+1},0)$
sending $\widehat M$ into itself such that ${\mathcal B}\circ
\widehat H=H\circ {\mathcal B}$ with $j_0^{\ell}\widehat
H=j_0^{\ell}{\rm Id}$.
\end{prop}

\begin{proof} Fix any integer $\ell\geq {\rm max}\, \{\alpha_n,3+6b_1\}$ and $H\colon (\C^{n+1},0)\to
(\C^{n+1},0)$  a holomorphic map sending $M$ into itself
satisfying $j_0^{\ell}H=j_0^{\ell}{\rm Id}$. We first show that
there exists a holomorphic map $\widehat H\colon (\C^{n+1},0)\to
(\C^{n+1},0)$ such that ${\mathcal B}\circ \widehat H=H\circ
{\mathcal B}$ with $j_0^{\ell}\widehat H=j_0^{\ell}{\rm Id}$. Then
we will show that the map $\widehat H$ is unique and sends
$\widehat M$ into itself.

We write the map $H=(F,G)\in \C^n\times \C$ and similarly the
desired map $\widehat H$ will be split as follows $(\widehat
F,\widehat G)$. We first note that since $(z,w)$ are normal
coordinates for $M$ and since $j_0^{\ell}H=j_0^{\ell}{\tt Id}$, we
have $G(z,w)=w(1+\Psi (z,w))$, where $\Psi$ is a convergent power series
vanishing at least to order $\ell -1$. Therefore setting
\begin{equation}\label{e:set}
\widehat G(z,w)=w\sqrt{1+(\Psi \circ {\mathcal
B})(z,w)}
\end{equation}
 we have that $\widehat G$ satisfies
\begin{equation}\label{e:firstgood}
\widehat G^2=G\circ {\mathcal B},\quad {\rm and}\quad \widehat
G(z,w)=w,\ {\rm up}\ {\rm to}\ {\rm order}\ \ell.
\end{equation}
Next, we claim that if we write $F=(F_1,\ldots,F_n)$, each power
series $F_j\circ {\mathcal B}$ is divisible by $w^{\alpha_{j}}$
for every $j=1,\ldots,n$. Indeed, by assumption we may write
$F_j(z,w)=z_j+\Phi_j(z,w)$, $\Phi_j$ being a power series
vanishing at least up to order $\ell$. Hence, since $(\Phi_j\circ
{\mathcal
B})(z,w)=\Phi_j(z_1w^{\alpha_1},\ldots,z_nw^{\alpha_n},w^2)$ and
since each $\alpha_j\geq 2$, we get that necessarily $w^{2(\ell
+1)}$ divides $(\Phi_j\circ {\mathcal B})(z,w)$ which proves the
claim since $2(\ell +1)>\alpha_n\geq \alpha_j$ for each $j$. We
may now set
\begin{equation}\label{e:cook}
\widehat F_j (z,w):=\frac{(F_j\circ {\mathcal B})(z,w)}{(\widehat G(z,w))^{\alpha_j}},
\end{equation}
which defines a convergent power series in view of \eqref{e:set}
and the above claim. From \eqref{e:firstgood} and \eqref{e:cook}
we automatically have $H\circ {\mathcal B}={\mathcal B}\circ
\widehat H$. To prove that for each $j=1,\ldots,n$, $F_j(z,w)=z_j$
up to order $\ell$, we write
$$\widehat F_j(z,w)=\frac{(F_j\circ {\mathcal B})(z,w)}{w^{\alpha_j}}\, \left(1+(\Psi\circ {\mathcal B})(z,w)\right)^{-\alpha_j/2},$$
and notice that $(F_j\circ {\mathcal B})(z,w)/w^{\alpha_j}$ agrees
with $z_j$ up to order (at least) $2(\ell+1)-1-\alpha_j\geq
2\ell+1-\alpha_n$ and that $\left(1+(\Psi\circ {\mathcal
B})(z,w)\right)^{-\alpha_j/2}$ agrees with $1$ up to order (at
least) $2\ell -1$. From this we conclude that $\widehat F_j(z,w)$
agrees with $z_j$ up to order at least $2\ell+1-\alpha_n\geq
\ell$, which proves that $j_0^{\ell}\widehat H=j_0^{\ell}{\rm
Id}$.

The uniqueness of $\widehat H$ is clear since there are only two
possibles choices for constructing a power series $T(z,w)$ such
that $T^{2}=G\circ {\mathcal B}$ i.e. $T=\pm \widehat G$ with
$\widehat  G$ as constructed above. Since we require $\widehat H$
to have the same $\ell$-jet as that of the identity mapping with
$\ell\geq \alpha_n\geq 1$, $T=\widehat  G$ is the only possible
choice for the normal component of $\widehat H$. This finally also
implies the uniqueness of the other components of $\widehat H$
since they must satisfy \eqref{e:cook}.

It remains to check that the constructed holomorphic map $\widehat H\colon (\C^{n+1},0)\to (\C^{n+1},0)$ indeed sends $\widehat M$ into itself. For this we first note that $M':=\widehat H(\widehat M)$ defines a  real-analytic hypersurface that is of infinite type since $\widehat M$
 is. Secondly, since $H$ sends $M$ into itself and ${\mathcal B}\circ \widehat H=H\circ {\mathcal
 B}$, we have that $\widehat H(\widehat M)\subset {\mathcal B}^{-1}(M)$. We now
 claim that $M'$ is a real-analytic hypersurface of the form ${\rm Im}\, w=({\rm Re}\, w)^{3+6b_1}\, \rho (z,\bar z,{\rm Re}\,
 w)$ for some convergent power series $\rho$. From the claim, one gets
 that $M'=\widehat M$ since $\widehat
 M$ is the unique real-analytic hypersurface of the above form
 contained in ${\mathcal B}^{-1}(M)$ by
 Proposition~\ref{p:construction}. Now, the proof of the claim is
 easily obtained by using the following fact that follows from our construction of the map $\widehat H$ : one has $\widehat G(z,w)=w+w^{2(\ell+1)}\delta
 (z,w)$ for some power series $\delta$. Since
 $2(\ell+1)>3+6b_1$, a direct computation of $\widehat H(\widehat
 M)$, left to the reader, shows that it has the form claimed above. The proof of Proposition~\ref{p:lift} is complete.
\end{proof}

Recall that $M$ (resp.\ $\widehat M$) has the finite jet
determination property at $0$ if there exists a positive integer
$r$ such that if $H\in \autM$ (resp.\ $H\in {\Aut (\widehat
M,0)}$) agrees with the identity mapping at the origin up to order
$r$, then $H$ is the identity. As an immediate consequence of
Proposition~\ref{p:lift}, we obtain the following crucial result.

\begin{cor}\label{c:endstory}
If $\widehat M$ has the finite jet determination property at $0$, so does $M$.
\end{cor}

\begin{rem}\label{r:soleremark}
If, throughout this section, $M$ is a formal real hypersurface (instead of a real-analytic one), still of the form
\eqref{e:normalformequation} with $R$ and each $\theta_j$ being formal power series, then Propositions \ref{p:construction} and \ref{p:lift} have a formal counterpart. The formal version of Proposition~\ref{p:construction} will provide a unique formal real hypersurface $\widehat M$ of the form \eqref{e:form1} with $\eta$ being a formal power series satisfying also \eqref{e:form2}. In addition, the lifting procedure given by Proposition~\ref{p:lift} also holds in that setting for formal holomorphic maps $H$.
The proofs of these statements are obtained by mimicing the proofs given here in a formal setting. The details are left to the reader.
\end{rem}

\subsection{Proof of Theorem \ref{t:main}}The proof of the theorem is achieved by combining
Proposition \ref{p:normalform}, Proposition \ref{p:construction},
Corollary \ref{c:endstory} and Proposition \ref{p:fjdgood}.


\bibliographystyle{plain}
\bibliography{bibliography_oct05}

\end{document}

%% file: definitions.tex
\newcommand{\inp}[2]{\langle #1 , #2 \rangle}

\newcommand{\C}{\mathbb{C}}

\newcommand{\R}{\mathbb{R}}

\newcommand{\N}{\mathbb{N}}

\newcommand{\dop}[1]{\frac{\partial}{\partial #1}}
\newcommand{\dopt}[2]{\frac{\partial #1}{\partial #2}}
\newcommand{\vardop}[3]{\frac{\partial^{|#3|} #1}{\partial {#2}^{#3}}}

\newcommand{\autM}{{\Aut (M,0)}}
\newcommand{\autMp}{{\Aut (M,p)}}

\newlength{\extendaxesby}

\DeclareMathOperator{\imag}{Im}
\DeclareMathOperator{\real}{Re}

\DeclareMathOperator{\Aut}{Aut}

\newtheorem{thm}{Theorem}
\newtheorem{lem}[thm]{Lemma}
\newtheorem{prop}[thm]{Proposition}
\newtheorem{cor}[thm]{Corollary}

\theoremstyle{definition}
\newtheorem{defin}[thm]{Definition}

\newtheorem{rem}[thm]{Remark}

\newtheorem{conj}[thm]{Conjecture}

\def \eps{\epsilon}